\title{Color-Critical Graphs Have Logarithmic Circumference}
\author{Asaf Shapira
\thanks{School of Mathematics and School of Computer Science, Georgia Institute of Technology, Atlanta, GA 30332. Supported in part by
NSF Grant DMS-0901355.}
\and Robin Thomas
\thanks{School of Mathematics, Georgia Institute
of Technology, Atlanta, GA 30332-0160.
Supported in part by NSF Grant number DMS-0739366.
21 August 2009.}
}

\date{}

\documentclass [letterpaper,11pt]{article}
\usepackage{amsfonts}
\usepackage{theorem}
\newtheorem{lemma}{Lemma}[section]
\newtheorem{coro}[lemma]{Corollary}
\newtheorem{theo}[lemma]{Theorem}

{\theorembodyfont{\rmfamily}
\newtheorem{definition}[lemma]{Definition}
}

\newcommand{\qed}{\hspace*{\fill} \rule{7pt}{7pt}}

\newcommand{\ignore}[1]{}

\def\mytextindent#1{\indent\indent\llap{\rm#1\enspace}\ignorespaces}
\def\myitem{\par\hangindent0pt\mytextindent}
\def\junk#1{}
\def\ell{l}

\begin{document}
\maketitle

\begin{abstract}
A graph $G$ is $k$-critical
if every proper subgraph of $G$ is $(k-1)$-colorable, but the graph
$G$ itself is not.
We prove that every $k$-critical graph on $n$ vertices has a cycle
of length at least $\log n/(100\log k)$,
improving a bound of Alon, Krivelevich and Seymour from 2000.
Examples of Gallai from 1963 show that the bound cannot be improved
to exceed $2(k-1)\log n/\log(k-2)$.
We thus settle the problem of bounding the minimal circumference of
$k$-critical graphs, raised by Dirac in 1952 and Kelly and Kelly in 1954.
\end{abstract}

\section{Introduction}\label{intro}

\junk{
A graph $G=(V,E)$ is said to be $k$-colorable if we can color the vertices of $G$ using the colors
$\{1,\ldots,k\}$ such that adjacent vertices receive distinct colors. The chromatic number of $G$, denoted $\chi(G)$, is the smallest
integer $k$ for which $G$ is $k$-colorable. The chromatic number of a graph is probably the most well studied graph invariant, see the book \cite{JT} for a comprehensive
discussion and references. As in many other areas, when studying a property of certain structures, it is natural and crucial to study the objects that are minimal with respect
to satisfying the property. When studying the chromatic number of a graph we define a graph to be $k$-critical if $\chi(G)=k$ and the removal of {\em any} edge from the graph results
in a new graph with chromatic number smaller than $k$. Clearly, every graph satisfying $\chi(G)=k$ contains as a subgraph a $k$-critical graph.
}

All graphs in this paper are finite and simple;
that is, they have no loops or multiple edges.
Paths and cycles have no ``repeated" vertices.
A graph $G$ is {\em $k$-critical}, where $k\ge 1$ is an integer,
if every proper subgraph of $G$ is $(k-1)$-colorable, but the graph
$G$ itself is not.
%We say that a graph is {\em critical} if it is $k$-critical for some $k\ge1$.
There is an easy description of $k$-critical graphs for $k\le3$, but
for $k\ge4$ their structure appears complicated and no meaningful
characterization is known.

The study of $k$-critical graphs was introduced in the $1940$s by Dirac
as part of his PhD Thesis.
Since then $k$-critical graphs have been studied
extensively, as documented for instance in~\cite[Chapter~5]{JT}.
In this paper we study the circumference of $k$-critical graphs,
where the circumference of a graph $G$ is the length of the longest
cycle in $G$.
The only $3$-critical graphs are odd cycles, but for $k\ge4$ the
circumference problem is more complicated.
For  integers $k\ge4$ and $n>k+1$ let $L_k(n)$ denote the largest integer~$l$
such that every $k$-critical graph on $n$ vertices has circumference
at least~$l$.
Elementary constructions show that the function $L_k(n)$ is well-defined
for all integers $k\ge4$ and $n>k+1$.
The study of the function $L_k(n)$ originated in the work of
Dirac~\cite{D1} and Kelly and Kelly~\cite{KK}.
%it is also formulated as~\cite[Problem~5.11]{JT}.

%of bounding $L_k(n)$ for $k \geq 4$ seems much more difficult, and
%from now on we always assume that $k \geq 4$ is a fixed integer and
%the number of vertices $n$ is much larger than $k$.
As every $k$-critical graph has minimum degree at least $k-1$, we
have $L_k(n) \geq k$. Dirac \cite{D1} showed that
$L_k(n) \geq 2k-2$ for all $n\ge 2k-2$ and conjectured that
$k$-critical graphs should contain much longer cycles.
Specifically, he conjectured that
for every fixed $k$ we have $\lim_{n \rightarrow \infty}L_k(n)= \infty$
and that actually $L_k(n) \geq c\sqrt{n}$.

The first non-trivial bounds on $L_k(n)$ where obtained in 1954 by
Kelly and Kelly \cite{KK} who showed that
%$\lim_{n\to\infty} L_k(n) =\infty$,
\begin{equation}\label{kelly}
\lim_{n\to\infty} L_k(n) =\infty\;,
\end{equation}
thus confirming the first conjecture of Dirac mentioned above.
According to~\cite{AKS} Kelly and Kelly~\cite{KK} actually proved that
\begin{equation}\label{kelly1}
L_k(n) \geq \sqrt{\frac{\log\log n}{\log\log\log n}}\;
\end{equation}
for every fixed $k\ge4$ and all sufficiently large $n$.
They also showed that
%$\liminf_{n\to\infty}L_4(n)/\log^2 n$
\begin{equation}\label{kelly2}
%\liminf_{n\to\infty}{L_4(n)\over\log^2 n}
\liminf_{n\to\infty}{L_4(n)/\log^2 n}\le 3/\log^2(27/4),
\end{equation}
thus disproving Dirac's second conjecture for $k=4$.
Dirac \cite{D} later extended the upper bound of \cite{KK} to
all~$k \geq 4$, and Read \cite{R} later improved the upper bound by
showing that
\begin{equation}
\label{read}
\liminf_{n\to\infty}L_k(n)/(\log n \cdot \log^{(2)} n \cdot \log^{(3)} n
\cdots \log^{(k-4)}n \cdot (\log^{(k-3)}n)^2)\le
(2/\log4)^{k-2},
\end{equation}
where $\log^{(i)}(x)$ is the $i$-times iterated logarithm function.
The best known upper bound on $L_k(n)$ was obtained in 1963
by Gallai \cite{G}, who improved and significantly simplified the previous
constructions by showing that for every integer $k\ge4$ there are
infinitely many integers $n$ such that

\begin{equation}\label{upperbound}
L_k(n) \leq \frac{2(k-1)}{\log (k-2)}\log n\;.
\end{equation}
We present Gallai's examples in Section~\ref{secgallai}, and we
also point out that the same graphs establish the related fact that
for every integer $k\ge4$ there are
infinitely many integers $n$ such that
\begin{equation}
\label{gallaipath}
\hbox{there exists a $k$-critical graph on $n$ vertices with no
{\em path} of length
exceeding } \frac{4(k-1)}{\log (k-2)}\log n\;.
\end{equation}

As for lower bounds on $L_k(n)$, the first improvement of the result of Kelly and Kelly \cite{KK} came after almost 50 years, when Alon, Krivelevich and Seymour \cite{AKS}
obtained the following (exponential) improvement of (\ref{kelly1})
for all integers $k\ge4$ and all integers $n\ge k+2$:
\begin{equation}\label{lowerbound}
L_k(n) \geq 2 \sqrt{\frac{\log (n-1)}{\log (k-2)}}\;.
\end{equation}

%Let us digress for a moment and discuss longest paths in $k$-critical graphs.
The proof of (\ref{lowerbound}) in~\cite{AKS} is based on a result (implicit in \cite{AKS})
which says that every $k$-critical graph on $n$ vertices has a path of length at least
$\log n/(\log(k-2))$. For completeness we state and prove it as Lemma~\ref{AKS} below.
This is asymptotically best possible by~(\ref{gallaipath}).

%\noindent
The main result of this paper is the following improvement of the
theorem of~\cite{AKS}.

\begin{theo}\label{theomain}
For every integer $k\ge4$ and every integer $n\ge k+2$ we have
$$
L_k(n) \geq \frac{\log n}{100 \log k}\;.
$$
\end{theo}
%Combining a variation on Gallai's upper bound (\ref{upperbound}) and our
%improved lower bound we get the following corollary, which
The following corollary solves,
for every fixed $k\ge4$, the problem of determining the order of
magnitude of $L_k(n)$.
The problem originated in the work of Dirac~\cite{D1} and
Kelly and Kelly~\cite{KK}, and is also stated in~\cite[Problem~5.11]{JT}.
%Dirac's problem and~\cite[Problem~5.11]{JT}.
The lower bound follows immediately from Theorem~\ref{theomain};
the upper bound follows by a minor modification of
Gallai's proof of~(\ref{upperbound}) and is presented
as Theorem~\ref{supergallai}.
%in Section~\ref{secgallai}.

\begin{coro}
\label{coro}
For every integer $k \geq 4$ and every integer $n\ge k+2$
$$
\frac{\log n}{100 \log k}\le
L_k(n)\le \frac{2(k-1)}{\log (k-2)}\log n+2k\;.
$$
\end{coro}

\noindent
%Dirac's problem asked for the circumference of $k$-critical
%graphs on $n$ vertices
%where $k$ is fixed and $n$ grows to infinity. This question
%is answered by the above corollary, but one can naturally ask
The corollary raises the obvious question
whether there exist
a function $f$ and absolute constants $c_1,c_2$ such that
$c_1f(k)\log n\le L_k(n)\le c_2f(k)\log n$.
%We do not know the answer.
This remains an interesting open problem.
The related (and perhaps easier) question, where we ask for the length
of the longest path, is also open.
Currently, the best known bounds for the latter problem are given
by~(\ref{gallaipath}) and Lemma~\ref{AKS}.

%In fact, it is open even for the related function $L'_k(n)$
%defined as the smallest
%integer $l$ such that every $k$-critical graph on $n$ vertices
%has a path on at least $l$ vertices.
%Our results imply that every bound for $L'_k(n)$ translates into
%a bound on $L_k(n)$.
%{\bf More details are called for.}

There is a related problem, formulated by Ne\v set\v ril and R\"odl
at the International Colloquium on Finite and Infinite Sets
in Keszthely, Hungary in 1973; see~\cite{JNbook} for a detailed
history of the problem.
Ne\v set\v ril and R\"odl asked whether it is true that for
every two integers $k,n\ge4$ there exists an integer $N$ such
that every $k$-critical graph on at least $N$ vertices has a
$(k-1)$-critical subgraph on at least $n$ vertices.
For $k=4$ the answer is yes, for the following reason.
By~(\ref{kelly}) a large enough $k$-critical graph $G$ has a long
cycle $C$. Since $G$ is not bipartite, it has an odd cycle, say $C'$.
The graph $G$ is $2$-connected by Lemma~\ref{2con}(i) below.
Now an elementary argument using just the $2$-connectivity of $G$
shows that $G$ has an odd cycle of length at least $|V(C)|/2$.
(The details may be found in~\cite{AKS}.)
This argument and Theorem~\ref{theomain} imply the following corollary.

\begin{coro}
\label{coro2}
Let $k,n\ge4$ be  integers. Then every $k$-critical graph on $n$
vertices has a $3$-critical subgraph on at least $\log n/(200\log k)$
vertices.
\end{coro}
This is an improvement over the bound $\sqrt{\log n/\log(k-1)}$
of Alon, Krivelevich and Seymour~\cite{AKS}.
The problem of Ne\v set\v ril and R\"odl is open for all $k\ge5$.

The rest of the paper is organized as follows.
In Section \ref{prelim} we prove a lemma that is implicit in~\cite{AKS}
and deduce~(\ref{lowerbound}) from it, and
give an overview of the proof of Theorem~\ref{theomain}.
In Section \ref{sec:lemmas}
we prove some basic results concerning $k$-critical graphs which
will be used in the proof of Theorem~\ref{theomain}.
%we prove some additional preliminary
%results concerning the structures we will be using in the main proof.
In Section~\ref{secbondy} we prove a variation of a theorem of
Bondy and Locke~\cite{BL}, stated as Theorem~\ref{BLtheo} below.
The proof of Theorem \ref{theomain} appears in Section~\ref{secmain}.
In Section~\ref{secgallai} we present Gallai's construction that leads
to the upper bound (\ref{upperbound}) and statement~(\ref{gallaipath}),
and point out how to deduce
the upper bound in Corollary~\ref{coro} from it.

\section{Proof Overview} % and Some Preliminary Results}
%%%%%%%%%%%%%%%%%%%%%%%%
\label{prelim}

%The main idea used in \cite{AKS} for proving of (\ref{lowerbound})
%is a variant of the following lemma, which we will also use later.

If $T$ is a tree and $x,y\in V(T)$, then there is a unique path
in $T$ with ends $x$ and $y$, and we will denote it by $xTy$.
Let $G$ be a graph, let $T$ be a spanning tree of $G$, and let
$v\in V(T)$. The tree $T$ is called a
{\em depth-first search (DFS) spanning tree rooted at $v$}
if for every edge $xy\in E(G)$ either $x\in V(vTy)$ or $y\in V(vTx)$.
It is easy to see that for every connected graph $G$ and every
vertex $v\in V(G)$ there is a DFS spanning tree in $G$ rooted at $v$.
For $X=\emptyset$ the following lemma is implicit in~\cite{AKS}.
The straightforward generalization will be needed later in the paper.
%We state and prove a variant of the lemma of \cite{AKS}
%since we will use it later in the paper.

\begin{lemma}\label{AKSmod}
Let $G$ be a $k$-critical graph on $n$ vertices, let $X\subseteq V(G)$
have size $s$, and let $T$ be a DFS spanning tree of $G\backslash X$
rooted at a vertex $t_0$.
Then for every integer $j\ge1$ the number of vertices of $T$ at distance
exactly $j$ from $t_0$ is at most $(k-1)^s$ if $j=1$ and
$(k-2)^{j-2}(k-1)^s$ otherwise.
%$j!(k-1)^x$ if $j\le k-3$ and at most
%$(k-2)!(k-2)^{j-k+2}(k-1)^x$ otherwise.
\end{lemma}

\paragraph{Proof.}
Let $G, X, T, t_0$ be as stated, let $X=\{x_1,x_2,\ldots,x_s\}$,
and let $t\in V(T)$ be a vertex at distance $j\ge1$ from $t_0$ in $T$.
We wish to define a sequence $Q(t)$.
Let $t_0,t_1,\ldots,t_{j-1},t_j=t$ be the vertices of the path $tTt_0$,
listed in order. Since $G$ is $k$-critical, the graph $G\backslash tt_{j-1}$
obtained from $G$ by deleting the edge $tt_{j-1}$ is $(k-1)$-colorable,
and since $t_0$ is adjacent to $t_1$ it has a $(k-1)$-coloring $\phi$ such that
$\phi(t_0)=1$ and $\phi(t_1)=2$.
We define $Q(t):=(\phi(t_2),\phi(t_3),\ldots,\phi(t_{j-1}),
\phi(x_1), \phi(x_2),\ldots,\phi(x_s))$.
Since for $i=1,2,\ldots,j$ the vertex $t_i$ is adjacent to $t_{i-1}$,
there are at most $(k-2)^{j-2}(k-1)^s$ sequences that arise this way
(at most $(k-1)^s$ if $j=1$).
%$\phi(t_i)$ can take on at most  $k-2$ different values.
It follows that if there are more than $(k-2)^{j-2}(k-1)^s$ vertices
at distance $j$ from $t_0$ in $T$ (or more than $(k-1)^s$ if $j=1$),
then there are two vertices $t,t'$ at distance $j$ from $t_0$ such that
$Q(t)=Q(t')$.
Let $t'_0=t_0,t'_1,\ldots,t'_{j-1},t'_j=t'$ be the vertices of the path
$t_0Tt'$,
let $p$ be the largest integer such that $t_p=t'_p$, and let the set
$Z$ consist of the vertex $t_{p+1}$ and
all its descendants in the rooted tree $(T,t_0)$.
Let $\phi$ be a coloring as above, and let $\phi'$ be the analogous
coloring of $G\backslash t't'_{j-1}$. The fact that $Q(t)=Q(t')$
implies that
\begin{equation}
\label{phiphi}
%$\phi(u)=\phi'(u)$ for every $u\in X\cup V(vTt_p)$.
\phi(u)=\phi'(u)\hbox{ for every }u\in X\cup V(t_0Tt_p).
\end{equation}
We now define a coloring $\psi$ by $\psi(u)=\phi(u)$ for every $u\in V(G)-Z$
and $\psi(u)=\phi'(u)$ for every $u\in Z$.
Since $T$ is a DFS spanning tree of $G\backslash X$, it follows that
every edge of $G$ with one end in $Z$ and the other end in $V(G)-Z$
has the other end in $X\cup V(t_0Tt_p)$. It follows from (\ref{phiphi})
that $\psi$ is a valid $(k-1)$-coloring of $G$, contrary to the
$k$-criticality of $G$.~\qed
\bigskip

The above lemma is the main tool in the proof of~(\ref{lowerbound}).
Alon, Krivelevich and Seymour~\cite{AKS} use it to deduce that every
$k$-critical graph has a long path, as follows.

\begin{lemma}
\label{AKS}
For every integer $k\ge4$ every $k$-critical graph on $n$ vertices has a
path of length at least $\log n /\log (k-2)$.
\end{lemma}

\paragraph{Proof.}
Let $k\ge4$ be an integer,  let $G$ be a $k$-critical graph on $n$
vertices, and let $v\in V(G)$ be some vertex in $G$.
Since $G$ is connected, it has
a DFS spanning tree $T$ rooted at $v$.
Let $h$ be the length of a longest path in $T$ with one end $v$.
By Lemma~\ref{AKSmod} applied with $X=\emptyset$ we deduce that
$$
n\le 1+1+1+(k-2)+(k-2)^2+\cdots+(k-2)^{h-2}\le\sum_{j=0}^{h-1}(k-2)^j\le
(k-2)^h,
$$
because $h\ge2$ and $k\ge4$.
It follows that $h\ge\log n /\log (k-2)$, as desired.~\qed
\bigskip

The next lemma is due to Dirac and Voss. A proof may be found
in~\cite{AKS,L}.

\begin{lemma}
\label{sqrtcycle}
If a $2$-connected graph has a path of length $l$, then it has
a cycle of length at least $2\sqrt{l}$.
\end{lemma}

Since every $k$-critical graph is $2$-connected by Lemma~\ref{2con} below,
the lower bound (\ref{lowerbound}) of
Alon, Krivelevich and Seymour
follows immediately from Lemmas~\ref{AKS} and~\ref{sqrtcycle}.
However, the bound in Lemma~\ref{sqrtcycle} is tight
and the bound in Lemma~\ref{AKS} is asymptotically tight by~(\ref{gallaipath}),
%\footnote{Although the bound in (\ref{upperbound}) was stated with relation
%to the circumference of the graph, Gallai's example actually gives a
%similar bound for the length of the longest path in a $k$-critical graph.
%See Section \ref{secgallai} for more details.}
and so an improvement of (\ref{lowerbound}) requires a different strategy.
For $3$-connected graphs the bound can be dramatically improved,
as shown by Bondy and Locke~\cite{BL}:

\begin{theo}
\label{BLtheo}
If a $3$-connected graph has a path of length $\ell$ then it has a
cycle of length at least $2\ell/5$.
\end{theo}

So combining Lemma~\ref{AKS} and Theorem~\ref{BLtheo} we get that every
$3$-connected $k$-critical graph has a cycle of length at least
$2\log n / (5\log k)$. Unfortunately, not all $k$-critical graphs
are $3$-connected, but those that are not can be constructed from
two smaller $k$-critical graphs. This is a result of Dirac~\cite{D2}
and is described in Lemma~\ref{2con}; a proof may also be found
in~\cite[Problem~9.22]{L}.
Thus one might hope that we could use this result of Dirac and
apply induction. That is indeed our strategy, but it turns out that
it is not enough (at least in our proof)
to just decompose the graph once and apply induction; instead,
we need to break the graph repeatedly by (non-crossing) cutsets of
size two and use the resulting tree-structure. That brings us to the
notion of tree-decomposition, which formalizes this break up process.

\begin{definition}
\label{deftreedecom}
A {\em tree decomposition} of a graph $G$ is a pair $(T,{\cal W})$
where $T$ is a tree and ${\cal W}=(W_t~:~t \in V(T))$ is a collection of
subsets of $V(G)$ such that
\begin{itemize}
\item $\bigcup_{t \in V(T)}W_t=V(G)$ and every edge of $G$ has both ends
in some $W_t$, and
\item If $t,t',t'' \in V(T)$ and $t'$ belongs to the unique path in $T$ connecting $t$ and $t''$, then $W_t \cap W_{t''} \subseteq W_{t'}$.
\end{itemize}
%We say that the tree-decomposition $(T,{\cal W})$ is {\em $2$-adhesive}
%if $|W_t\cap W_{t'}|=2$ for every edge $tt'\in E(T)$.
For $t\in V(T)$ we define the {\em torso of $(G,T,{\cal W})$ at $t$}
to be the graph with vertex-set $W_t$ in which $u,v\in W_t$ are adjacent
if either they are adjacent in $G$ or $u,v\in W_{t'}$ for some
neighbor $t'$ of $t$ in $T$.
We say that the tree-decomposition $(T,{\cal W})$ is {\em standard}
if $|W_t\cap W_{t'}|=2$ for every edge $tt'\in E(T)$
and each torso of $(G,T,{\cal W})$ is $3$-connected or a cycle.
(A graph $G$ is $t$-connected if it has at least $t+1$ vertices
and $G\backslash X$ is connected for every set $X\subseteq V(G)$ of
size at most $t-1$).
\end{definition}

We will see in Lemma~\ref{standard} that every $k$-critical graph
has a standard tree-decomposition.
% such that every torso is $3$-connected or a cycle.
The torsos are not necessarily critical, but they are very close,
so that is not really an issue, and so for the purpose of this outline we can
pretend that every torso satisfies the conclusion of Lemma~\ref{AKS}.
By Theorem~\ref{BLtheo} we deduce that each torso has a sufficiently
long cycle, but we need more.
We actually need a ``linkage", a set of two disjoint paths with prescribed
ends so that we can combine these linkages in individual torsos
to produce a cycle in the original graph.
We deduce the existence of such a linkage of desired length from
Theorem~\ref{BLtheo} in Lemma~\ref{disjointpaths}, but only under the
assumption that the
two sets of prescribed ends are disjoint from each other;
otherwise Lemma~\ref{disjointpaths} is false.
Thus when the two sets of prescribed ends are not disjoint
we need a different method.
In that case we are really looking for one path rather than
a linkage, and we
use Lemma \ref{AKSmod} to find it.

%use Lemmas~\ref{longpath} and~\ref{sptree} instead.

\section{Some Preliminary Lemmas}
%%%%%%%%%%%%%%%
\label{sec:lemmas}

In this section we prove some preliminary lemmas and introduce several concepts that will be used later
on in the proof of Theorem \ref{theomain} in Section \ref{secmain}. Our first lemma is well-known and appears in~\cite[Exercise~12.20]{Die}.
%We say that a tree-decomposition $(T,{\cal W})$ of a graph $G$ is
%{\em standard} if it is $2$-adhesive and every torso of $(T,{\cal W})$
%is either $3$-connected or a cycle.

\begin{lemma}
\label{standard}
Every $2$-connected graph has a standard tree-decomposition.
\end{lemma}

Let $G$ be a graph, and let $x,y$ be distinct vertices of $G$.
If $x,y$ are not adjacent, then we define $G+xy$ to be the graph
obtained from $G$ by adding an edge joining $x$ and $y$, and if
$x,y$ are adjacent, then we define $G+xy$ to be $G$.
We define $G/xy$ to be the graph obtained from $G$ by deleting the
edge $xy$ (if it exists), identifying the vertices $x$ and $y$ and
deleting all resulting parallel edges.
Actually, in all applications of the operation $G/xy$ in this paper
the vertices $x$ and $y$ will not be adjacent, and they will have
no common neighbors, so the clause about deleting parallel edges
will not be needed.
Statement (iv) of the following lemma is due to Dirac~\cite{D2},
and a proof may also be found in~\cite[Problem~9.22]{L}.

\begin{lemma}
\label{2con}
Let $k\ge4$ be an integer,  let $G$ be a $k$-critical graph, and
let $u,v\in V(G)$ be such that $G\backslash\{u,v\}$, the
graph obtained from $G$ by deleting the vertices $u$ and $v$, is disconnected.
Then
\myitem{(i)} $u\ne v$, and hence $G$ is $2$-connected,
\myitem{(ii)} $u$ is not adjacent to $v$,
\myitem{(iii)} $G\backslash\{u,v\}$ has exactly two components, and
\myitem{(iv)} there are unique proper induced subgraphs $G_1,G_2$ of
$G$ such that
$G=G_1 \cup G_2$, $V(G_1)\cap V(G_2)=\{u,v\}$, the graphs
$G_1\backslash\{u,v\}$ and $G_2\backslash\{u,v\}$ are the two components
of $G\backslash\{u,v\}$, $u,v$ have no common neighbor in $G_2$,
and $G_1+uv$ and $G_2/uv$ are $k$-critical.
\end{lemma}

\paragraph{Proof.}
Since $G\backslash\{u,v\}$ is disconnected, there is an integer $s\ge2$
such that the graph $G$ can be expressed as $G=G_1\cup G_2\cup\cdots\cup G_s$,
where the graphs $G_i$ are pairwise edge-disjoint, each has at least
one edge, and $V(G_i\cap G_j)=\{u,v\}$
for distinct integers $i,j\in\{1,2,\ldots,s\}$.
Since $G$ is $k$-critical, each $G_i$ is $(k-1)$-colorable.
If $u=v$, then each $G_i$ has a $(k-1)$-coloring that gives the
vertex $u$ color $1$. Those colorings can be combined to produce
a $(k-1)$-coloring of $G$, contrary to the $k$-criticality of $G$.
Thus $u\ne v$, and statement (i) follows.

%If $u$ and $v$ are adjacent, then each $G_i$ has a $(k-1)$-coloring that
%gives the vertices $u,v$ the colors $1$ and $2$, respectively.
%Those colorings can be combined to produce a $(k-1)$-coloring of $G$
%as above, a contradiction. This proves (ii).

We now prove (ii) and (iii) simultaneously. If one of them does not
hold, then we may assume that $s=3$.
(If (ii) does not hold, then $G_3$ may be chosen to consist
of $u$, $v$, and the edge joining them.)
Let us say that $G_i$ is of type one if some $(k-1)$-coloring of $G_i$
gives $u$ and $v$ the same color, and let us say it is of type two
if some $(k-1)$-coloring of $G_i$
gives $u$ and $v$ different colors.
Thus each $G_i$ is of type one or type two. We claim that no $G_i$
is of both types.
For suppose for a contradiction that $G_3$ is of both types.
But $G_1$ and $G_2$ are of the same type, because $G_1\cup G_2$
is $(k-1)$-colorable by the $k$-criticality of $G$, and hence
it follows that $G_1\cup G_2\cup G_3=G$ is $(k-1)$-colorable,
a contradiction.
Thus no $G_i$ is of both types.
We may assume that $G_1$ and $G_2$ are of the same type and
that $G_3$ is of different type. But then $G_1\cup G_3$ is not
$(k-1)$-colorable, contrary to the $k$-criticality of $G$.
This proves (ii) and (iii). In particular, $s=2$.

We now prove (iv).
%The uniqueness provision follows from (iii), and so it suffices to
%prove that $G_1$ and $G_2$ exist, $u,v$ have no common neighbor in $G_2$,
%and $G_1+uv$ and $G_2/uv$ are $k$-critical.
We may assume from the symmetry that $G_1$ is of type one and $G_2$
is of type two.
Since $G$ is not $k$-colorable, it follows that $G_1$ is not of
type two and $G_1$ is not of type one.
%Since $G$ is not $(k-1)$-colorable, we may assume from the symmetry
%that every $(k-1)$-coloring of $G_1$ gives $u$ and $v$ the same
%color, and that every $(k-1)$-coloring of $G_1$ gives $u$ and $v$
%distinct colors. (Otherwise some $(k-1)$-colorings of $G_1$ and
%$G_2$ may be combined to produce a $(k-1)$-coloring of $G$ as above.)
It follows that neither $G_1+uv$ nor $G_2/uv$ is $(k-1)$-colorable.
It remains to show that $u,v$ have no common neighbor in $G_2$
and that every proper subgraph of $G_1+uv$ and $G_2/uv$ is $(k-1)$-colorable,
because the remaining properties are clear or follow from (iii).
Suppose for a contradiction that $w\in V(G_2)-\{u,v\}$ is a neighbor
of both $u$ and $v$ in $G$.
Then the graph $G\backslash uw$ has a $(k-1)$-coloring $\phi$ by the
$k$-criticality of $G$. Since $G_1$ is not of type two we deduce
that $\phi(u)=\phi(v)$. But $\phi(v)\ne\phi(w)$, because $v$ is adjacent
to $w$ in $G\backslash uw$. Thus $\phi$ is $(k-1)$-coloring of $G$,
a contradiction.
This proves that $u,v$ have no common neighbor in $G_2$.
Let $e$ be an edge of $G_1+uv$.
If $e\ne uv$, then let $\psi$ be a $(k-1)$-coloring of $G\backslash e$.
We have $\psi(u)\ne\psi(v)$, because $G_2$ is not of type one,
and hence $\psi$ is a $(k-1)$-coloring of $(G_1+uv)\backslash e$.
For $e=uv$ we note that $G_1$ is $(k-1)$-colorable by the $k$-criticality
of $G$.
Finally, let $f$ be an edge of $G_2/uv$, and let $\lambda$ be a
$(k-1)$-coloring of $G\backslash f$. Since $G_1$ is not of type two,
$\lambda(u)=\lambda(v)$, and hence $\lambda$ can be converted to
a $(k-1)$-coloring of $(G_2/uv)\backslash f$, as desired.~\qed
\bigskip

The first three statements of Lemma~\ref{2con} have the following
consequence.

\begin{lemma}
\label{virtue}
Let $k\ge4$ be an integer, let $G$ be a $k$-critical graph, and let
$(T,{\cal W})$ be a standard tree-decomposition of $G$.
% such that each $W_t$ has at least three elements.\marginpar{or standard?}
\myitem{(i)} If $t_0,t_1$ are adjacent in $T$, then the two
vertices in the set $W_{t_0}\cap W_{t_1}$ are not adjacent, and
\myitem{(ii)}
if $t_1,t_2$ are distinct neighbors of $t_0$ in $T$, then
$W_{t_0}\cap W_{t_1}\ne W_{t_0}\cap W_{t_2}$.
\end{lemma}

\paragraph{Proof.}
We prove only (ii), leaving (i) to the reader.
Suppose for a contradiction that $W_{t_0}\cap W_{t_1}=W_{t_0}\cap W_{t_2}$,
and let $X$ denote this $2$-element set.
Let $i\in\{0,1,2\}$. Since $(T,{\cal W})$ is standard, $W_{t_i}$ has
at least three elements, and hence
%By hypothesis
there exists a vertex $v_i\in W_{t_i}-X$.
Then $v_0,v_1,v_2$ belong to three different components of $G\backslash X$,
contrary to Lemma~\ref{2con}(iii).~\qed
\bigskip

Part (iv) of Lemma~\ref{2con} leads to the following construction,
which will modify each torso of a tree-decomposition of a $k$-critical
graph and turn it into a $k$-critical graph.
Let $G$ be a $k$-critical graph and let  $(T,{\cal W})$ be a
standard tree-decomposition of $G$ such that
each $W_t$ has at least three elements.
Let $t\in V(T)$, and let
$u,v\in W_t$ be distinct.
We say that the pair $uv$ is a {\rm virtual edge} of $W_t$ if
$W_t\cap W_{t'}=\{u,v\}$ for some neighbor $t'$ of $t$ in $T$.
Thus Lemma~\ref{virtue} asserts that the virtual edges of each $W_t$
are pairwise distinct, and that they are not edges of $G$
(but they are edges of the torso at $t$, by definition of torso).
We now classify virtual edges of $W_t$ into additive and contractive,
as follows.
Let $uv$ be a virtual edge of $W_t$, and let $t'$ be the neighbor of $t$
in $T$ such that $W_t\cap W_{t'}=\{u,v\}$.
Since $G\backslash\{u,v\}$ is disconnected, there exist graphs
$G_1,G_2$ as in Lemma~\ref{2con}(iv).
Then $W_t$ is a subset of exactly one of $V(G_1), V(G_2)$;
if $W_t\subseteq V(G_1)$, then we say that the virtual edge $uv$ is
{\em additive}; otherwise we say that it is {\em contractive}.
We now define a graph $N_t$ as the graph
obtained from $G[W_t]$ by adding the edge $uv$ for every additive virtual
edge $uv$ of $W_t$, and identifying the vertices $u$ and $v$ for every
contractive virtual edge $uv$ of $W_t$.
In other words, $N_t$ can be regarded as being obtained from the torso
of $(G,T,{\cal W})$ at $t$ by contracting all
contractive virtual edges of $W_t$.
We call $N_t$ the {\em nucleus} of $(G,T,{\cal W})$ at $t$.
The next lemma shows that the nucleus is well-defined in the sense
that the vertex identifications used during the construction do not
produce loops or parallel edges.

\junk{
We denote by $N_t$ the {\em nucleus} at a vertex $t \in T$ and define
it inductively as follows;
if $|T|=1$ then $N_t$ is $G$ in which case we say that $N_t$ is obtained
from $W_t$ in 0 steps. Otherwise, pick any neighbor $y$ of $t$ in $T$,
let $\{x,y\}=W_t \cap W_y$ and let $T_t$ and $T_y$ be the two trees
resulting from the removal of the edge $(t,y)$ from $T$, in such a way that $t \in T_t$ and $y \in T_y$. Finally, set $G_t=G[\bigcup_{v \in T_t} W_v]$ and $G_y=G[\bigcup_{v \in T_y} W_v]$. By the third assertion of Lemma \ref{2con}, either $G'=G_t+(x,y)$ or $G'=G_t/\{x,y\}$ is $k$-critical. Let $(T',{\cal W'})$ be a tree-decomposition of $G'$ which is obtained from $(T,W)$ (the tree decomposition of $G$) as follows: $T'$ is the tree $T_t$ defined
above and ${\cal W'}=\{W'_v~:~ v \in T_t\}$ where $W'_v=W_v$ for every $v \neq t$ and $W'_t=W_t$ if $G_t+(x,y)$ is a $k$-critical graphs, while
$W'_t$ is obtained from $W_t$ by merging together $\{x,y\}$ if $G_t/\{x,y\}$ is $k$-critical. By the third assertion of Lemma \ref{2con}, the graph $G'$ is $k$-critical. We thus define the nucleus of $W_t$ as the nucleus of $W'_t$. In this case we say that $N_t$ is obtained from
$W_t$ by $k+1$ steps, where $k$ is the number of steps needed to get $N_t$ from $W'_t$.
}

\begin{lemma}
\label{nucleus}
Let $k\ge4$ be an integer, let $G$ be a $k$-critical graph, let
$(T,{\cal W})$ be a standard tree-decomposition of $G$,
let $t\in V(T)$ and let $H$ denote the torso of $(G,T,{\cal W})$ at $t$.
Then
\myitem{(i)} the subgraph of $H$ induced
by contractive virtual edges of $W_t$ is a forest, and for every component
$R$ of this forest and every $v\in V(H)-V(R)$, at most one vertex of
$R$ is adjacent to $v$ in $H$, and
\myitem{(ii)}
the nucleus $N$ of $(G,T,{\cal W})$ at $t$ is $k$-critical.
\end{lemma}

\paragraph{Proof.}
We proceed by induction on the number of vertices of $T$.
If $T$ has only one vertex, then there are no virtual edges and $N_t=G$,
and hence both statements of the lemma hold.
We may therefore assume that $T$ has more than one vertex, and that
the lemma holds for all $k$-critical graphs that have a
standard tree-decomposition
using a tree with strictly fewer than $|V(T)|$ vertices.
Let $t'$ be a neighbor of $t$ in $T$, and let $W_t\cap W_{t'}=\{u,v\}$,
so that $uv$ is a virtual edge of $W_t$.
Let $T'$ be the component of $T\backslash tt'$ containing $t$,
let ${\cal W}'=(W_r: r\in V(T'))$, and let $G'$ be the subgraph of $G$
induced by the union of all $W_r$ over all $r\in V(T')$.

Assume first that $uv$ is an additive virtual edge.
Then $G'+uv$ is $k$-critical by Lemma~\ref{2con}(iv)
and $(T',{\cal W}')$ is a standard tree-decomposition of $G'+uv$, where
$T'$ has strictly fewer vertices than $T$.
Furthermore, $H$ is equal to the
torso of $(G'+uv,T',{\cal W}')$ at $t$,
and $N$ is equal to the nucleus of $(G'+uv,T',{\cal W}')$ at $t$.
Thus both conclusions follow by induction applied to $G'+uv$ and
the tree-decomposition $(T',{\cal W}')$.
This completes the case when $uv$ is an additive virtual edge.

We may therefore assume that $uv$ is a contractive virtual edge.
In this case we proceed analogously, applying induction to the
graph $G'/uv$ and the tree-decomposition obtained from $(T',{\cal W}')$
by replacing each occurrence of $u$ or $v$ by the new vertex of $G'/uv$
that resulted from the identification of $u$ and $v$.
In the proof of (i) we take advantage of the provision in
Lemma~\ref{2con}(iv) that guarantees that $u,v$ have no common neighbor
in $G'$.~$\qed$
\bigskip

\begin{lemma}
\label{ENEH}
If a graph $N$ is obtained from a graph $H$ by repeatedly contracting edges,
each time contracting an edge that belongs to no triangle, and $H$ has
minimum degree at least three, then $|E(H)|\le3|E(N)|$.
\end{lemma}

\paragraph{Proof.}
Let $d_1,d_2\ldots,d_n$ be the degree sequence of $N$, and let us
consider the reverse process that produces $H$ starting from $N$.
Then the $i^{th}$ vertex of $N$ gives rise to at most $d_i-3$
new edges of $H$. Thus
$$
|E(H)| \le |E(N)| + (d_1-3) + (d_2-3) +...+ (d_n-3) \le 3|E(N)|.
$$
as desired.~\qed

\begin{lemma}
\label{3connbags}
Let $k\ge4$ be an integer, let $G$ be a $k$-critical graph, let
$(T,{\cal W})$ be a standard tree-decomposition of $G$, and let
$t\in V(T)$.
Then the torso of $(G,T,{\cal W})$ at $t$ is $3$-connected.
\end{lemma}

\paragraph{Proof.}
Let $H$ denote the torso of $(G,T,{\cal W})$ at $t$.
% and let $N$ be the nucleus of $(G,T,{\cal W})$ at $t$.
If $H$ is not $3$-connected, then it is a cycle by the definition
of standard tree-decomposition. But then
the nucleus of $(G,T,{\cal W})$ at $t$ is a cycle, because it
is obtained from $H$ by contracting edges, contrary to
Lemma~\ref{nucleus}(ii).~\qed
\bigskip

\begin{lemma}
\label{degreeT}
Let $k\ge4$ be an integer,
let $G$ be a $k$-critical graph,
let $(T,{\cal W})$ be a standard tree-decomposition of $G$,
let $t\in V(T)$, and let $N$ be the nucleus of $(G,T,{\cal W})$ at $t$.
Then
$
\deg_T(t) \leq 3|E(N)|
$.
%where $deg(t)$ is the degree of $t$ in $T$, and
%where $N_t$ is the nucleus at $t$.
\end{lemma}

\paragraph{Proof.}
Let $H$ be the torso of $(G,T,{\cal W})$ at $t$.
We first notice that $\deg_T(t) \leq |E(H)|$, because each neighbor
of $t$ in $T$ gives rise to a unique virtual edge of $G$ at $t$
by Lemma~\ref{virtue}(ii), and each virtual edge belongs to $H$.
The graph $H$ is $3$-connected by Lemma~\ref{3connbags}.
By Lemma~\ref{nucleus}(i) the graph $N$ is obtained
from $H$ as in Lemma~\ref{ENEH},
and hence $|E(H)|\le3|E(N)|$ by that lemma, as desired.~\qed

\begin{lemma}
\label{edgesinnuc}
Let $k\ge4$ be an integer,
let $G$ be a $k$-critical graph,
let $(T,{\cal W})$ be a standard tree-decomposition of $G$, and
let $t\in V(T)$.
% and let $N$ be the nucleus of $(G,T,{\cal W})$ at $t$.
Then the nucleus of $(G,T,{\cal W})$ at $t$
has at least as many edges as $G[W_t]$.
\end{lemma}

\paragraph{Proof.}
This follows from the fact that no edges of $G[W_t]$ are lost during the
construction of the nucleus.~\qed
\bigskip

\begin{lemma}
\label{pathtorso}
Let $k\ge4$ be an integer,
let $G$ be a $k$-critical graph,
let $(T,{\cal W})$ be a standard tree-decomposition of $G$,
let $t\in V(T)$, and
let $N$ be the nucleus of $(G,T,{\cal W})$ at $t$.
Then the torso of $(G,T,{\cal W})$ at $t$
has a path of length at least $\frac12\log |E(N)|/\log k$.
\end{lemma}

\paragraph{Proof:}
Let $H$ denote the torso of $(G,T,{\cal W})$ at $t$.
By Lemma~\ref{nucleus} the graph $N$ is $k$-critical,
and so by Lemma~\ref{AKS} it has a path of length at least
$\log |V(N)|/\log k\ge\frac{1}{2}\log |E(N)|/\log k$.
Since $N$ is obtained from $H$ by contracting edges, $H$ has a path
at least as long.~$\qed$
\bigskip

Finally we need an easy lemma about trees.
If $G$ is a graph, $\phi: V(G) \mapsto \{0,1,2,3,\ldots\}$ is a mapping,
and $H$ is a subgraph of $G$, then we define
$\phi(H):=\sum_{v\in V(H)} \phi(v)$.

\begin{lemma}
\label{longpath}
Let $k\ge2$ be an integer, let $T$ be a tree, let $r\in V(T)$, and assume
that for every integer $\ell \geq 1$ there are at most $k^{\ell}$
vertices at distance exactly $\ell$ from $r$ in $T$.
Let $\phi: V(T) \mapsto \{0,1,\ldots\}$ be a weight function with
$\phi(r)=0$ and $\phi(t) \neq 0$ for at least one vertex in $t\in V(T)$.
Then there
exists a vertex $t\in V(T)$ at distance exactly $\ell$ from $r$ in $T$
such that $\phi(t)>0$ and
$$
2\ell\log k + \log \phi(t) \geq \log \phi(T)\;.
$$
\end{lemma}

\paragraph{Proof.}
For every integer $\ell \geq 0$
let $D_{\ell}$ be the set of vertices of $T$ at distance
exactly $\ell$ from $r$.
Since $\phi(r)=0$ we have that for some $\ell \geq 1$
$$
\sum_{t \in D_{\ell}} \phi(t) \geq \phi(T)/2^{\ell}\;,
$$
since if this is not the case, then
$$
\phi(T)=\sum_{\ell \geq 1}\sum_{t \in D_{\ell}} \phi(t) < \phi(T)
\sum_{\ell \geq 1}2^{-\ell} \leq \phi(T)\;,
$$
a contradiction.
Since $|D_{\ell}| \leq k^{\ell}$, we deduce that there is a vertex
$t\in D_\ell$ satisfying
\begin{equation}\label{depth}
\phi(t) \geq \phi(T)/2^{\ell}k^{\ell} \geq \phi(T)/k^{2\ell}\;,
\end{equation}
because $k\ge2$.
It follows that
$$
2\ell\log k + \log \phi(t)\geq
2\ell\log k + \log \phi(T) - 2\ell\log k=\log \phi(T)\;,
$$
as desired.~$\qed$

\section{An Application of the Theorem of Bondy and Locke}
%%%%%%%%%%%%%%%%%%%%%%%%%%%%%%%%%%%%%%%%%%%%%%%%%%%%%%%%%
\label{secbondy}

Let $G$ be a graph, and let $X,Y\subseteq V(G)$ be disjoint sets of
size two. A {\em linkage} in $G$ from $X$ to $Y$ is a set $\{P_1,P_2\}$
of two disjoint paths, each with one end in $X$ and the other end in $Y$.
The {\em length} of the linkage is defined to be $|E(P_1)|+|E(P_2)|$.
The following is the main result of this section.

\begin{lemma}
\label{disjointpaths}
Let $G$ be a $3$-connected graph, let
$X,Y\subseteq V(G)$ be disjoint sets of
size two, and suppose that
$G$ has a path of length at least $\ell$. Then $G$ has
a linkage from $X$ to $Y$ of length at least $\ell/25$.
\end{lemma}

The assumption that the sets $X,Y$ be disjoint is necessary, as the
following example shows.
Let $t\ge2$ be an integer, and let $H$ be
the graph obtained from a path with vertex-set $\{v_0,v_1,\ldots,v_{t^2}\}$
in order by adding an edge joining $v_{it}$ and $v_{(i+2)t}$ for
every $i=0,1,\ldots,t-2$.
Finally, let $G$ be obtained from $H$ by adding a vertex $u$ joined
to all vertices of $H$.
One can easily verify that $G$ is $3$-connected, has a path of length $t^2$, and yet
every linkage from $\{u,v_0\}$ to $\{u,v_2\}$ has length at most linear in $t$.
%the graph obtained from disjoint paths $P_0,P_1,\ldots,P_t$, each on
%$t$ vertices, by identifying the $i^{th}$ vertex of $P_0$ with
%one end of $P_i$.
%Finally, let $G$ be obtained from $H$ by adding three vertices $x,y,z$
%adjacent to each other and to every vertex of $H$.
%Thus $G$ is clearly $3$-connected.
%If $X=\{x,z\}$ and $Y=\{y,z\}$, then every linkage from $X$ to $Y$ in $G$
%has length at most $3t-1$, and yet $G$ has $t^2+3$ vertices.

In the proof of Lemma~\ref{disjointpaths}
we will make use of the following lemma,
which follows from the standard ``augmenting path" proof of
Menger's theorem or the Max-Flow Min-Cut Theorem;
see, for instance~\cite[Section~3.3]{Die}.

\begin{lemma}
\label{augment}
Let $r\ge1$ be an integer, let $G$ be an $r$-connected graph,
let $S$ and $T$ be two subsets of the vertex-set of $G$, each of size
at least $r$,
and let $P_1,P_2,\ldots,P_{r-1}$ be disjoint paths
such that for $i=1,2,...,r-1$, the path $P_i$
has ends $s_i \in S$ and $t_i \in T$.
%If there exist $r+1$ disjoint paths between $S$ and $T$ in $G$,
Then there exist disjoint paths $Q_1,Q_2,\ldots,Q_{r}$ in $G$
between $S$ and $T$ in such a way that all but one of
the paths $Q_i$ has an end in $\{s_1,s_2,...,s_{r-1}\}$,
and all but one of the paths
$Q_i$ has an end in $\{t_1,t_2,...,t_{r-1}\}$.
\end{lemma}

The proof of Lemma \ref{disjointpaths} will consist of three steps.
In the first step we will obtain either
a required  linkage, or a similar  structure we call {\em hammock},
which we introduce next.
In the second step we show that if a $3$-connected graph has a long hammock,
then it has either a long linkage, or a long ``non-singular" hammock.
Finally, we  show how to get a
required long linkage from the existence of a long non-singular hammock.

A {\em hammock} in $G$ from $X$ to $Y$ is a quadruple
$\eta=(P_1,P_2,R_1,R_2)$, where
\myitem{$\bullet$} $\{P_1,P_2\}$ is a linkage from $X$ to $Y$, where
$P_i$ has ends $x_i\in X$ and $y_i\in Y$,
\myitem{$\bullet$} $R_i$ is a path with ends $s_i\in V(P_1)$
and $t_i\in V(P_2)$, and is otherwise disjoint from $P_1\cup P_2$,
\myitem{$\bullet$} the paths $R_1,R_2$ are disjoint, except possibly
$s_1=s_2$,
\myitem{$\bullet$} the vertices $x_1,s_1,s_2,y_1$ occur on $P_1$ in
the order listed (but are not necessarily distinct), and
\myitem{$\bullet$} the vertices $x_2,t_1,t_2,y_2$ occur on $P_2$ in
the order listed (but are not necessarily distinct).

\noindent The {\em length} of the hammock $\eta$ is
defined to be $|E(R_1)|$.
%(This may seem strange, but it seems to work best for our purposes.)
(It may seem more natural to define the length of $\eta$ to be
$|E(R_1)|+|E(R_2)|$. Indeed, by doing so it is possible to improve the
constant $25$ in Lemma~\ref{disjointpaths} to $17.5$, but only at the
expense of more extensive case analysis.
The extra effort did not seem justified.)
We say that $\eta$ is {\em singular} if $s_1=s_2$, and {\em non-singular}
otherwise.

Let us recall that if $P$ is a path and $u,v\in V(P)$,
then by $uPv$ we denote the
unique subpath of $P$ with ends $u$ and~$v$.

\begin{lemma}
\label{lem:A}
Let $G$ be a $3$-connected graph, let $X,Y\subseteq V(G)$ be
disjoint sets of size two, and assume that $G$ has a cycle $C$ of
length $l$. Then $G$ has either a linkage from $X$ to $Y$
of length  at least $l/5$,
or a hammock from $X$ to $Y$
or from $Y$ to $X$ of length at least $2l/5$.
\end{lemma}

\paragraph{Proof:}
Assume first that there exist four disjoint paths $P_1,P_2,P_3,P_4$,
each with one end in $X\cup Y$ and the other end in $V(C)$.
For $i=1,2,3,4$ let $u_i$ and $v_i$ be the ends of $P_i$
such that $u_i\in V(C)$, $v_1,v_2\in X$ and $v_3,v_4\in Y$.
If $u_1,u_3,u_2,u_4$ occur on $C$ in the order listed,
then let $C_{13}$ denote the subpath of $C\backslash\{u_2,u_4\}$ with
ends $u_1$ and $u_3$, and let $C_{23}$, $C_{24}$, and $C_{14}$
be defined analogously.
Then either $P_1\cup P_2\cup P_3\cup P_4\cup C_{13}\cup C_{24}$
or $P_1\cup P_2\cup P_3\cup P_4\cup C_{23}\cup C_{14}$
is a linkage from $X$ to $Y$ of length at least $\ell/2$, as desired.
Thus we may assume that $u_1,u_2,u_3,u_4$ occur on $C$ in the order listed.
Using analogous notation,
if $|E(C_{14})|+|E(C_{23})|\ge \ell/5$, then
$P_1\cup P_2\cup P_3\cup P_4\cup C_{14}\cup C_{23}$ is a linkage
from $X$ to $Y$ in $G$ of length at least $\ell/5$.
Thus we may assume that $|E(C_{12})|+|E(C_{34})|\ge 4\ell/5$, and
so from the symmetry we may assume that $|E(C_{12})|\ge 2\ell/5$.
Then $(P_1\cup C_{14}\cup P_4,P_2\cup C_{23}\cup P_3,C_{12},C_{34})$
is a hammock from $X$ to $Y$ in $G$ of length at least $2\ell/5$, as desired.
This completes the case when $G$ has four disjoint paths from $X\cup Y$
to $V(C)$.

We may therefore assume that those four paths do not exist,
and hence by Menger's theorem $G$ can be expressed as $G_1\cup G_2$,
where $|V(G_1)\cap V(G_2)|=3$,
$X\cup Y\subseteq V(G_1)$ and $V(C)\subseteq V(G_2)$.
%We may assume that subject to the conditions listed, $G_1$ and $G_2$
%are chosen so that $G_1$ is minimal with respect to taking subgraphs.
Since $G$ is $3$-connected there exist three disjoint paths $P_1,P_2,P_3$
from $X\cup Y$ to $V(C)$ with no internal vertices in $X\cup Y\cup V(C)$.
By symmetry, we may assume that
$P_i$ has ends $u_i$ and $v_i$, where $u_i\in V(C)$, $v_1,v_2\in X$
and $v_3\in Y$.
Then for $i=1,2,3$ the set $V(G_1)\cap V(G_2)\cap V(P_i)$ includes
a unique vertex, say $w_i$.
Thus $V(G_1)\cap V(G_2)=\{w_1,w_2,w_3\}$.
Please note that the sets $\{w_1,w_2,w_3\}$ and
$\{u_1,u_2,u_3,v_1,v_2,v_3\}$ may intersect.

%Since a subpath of $P_3$ joins $v_3$ to $w_3$ in $G_1\backslash V(P_1\cup P_2)$,
%the minimality of $G_1$ and
%the ``augmenting path" proof of Menger's
%theorem or the max-flow min-cut theorem
%(see, for instance,~\cite[Section~3.3]{Die})
%applied to augment the path $v_3Pw_3$

By Lemma~\ref{augment} applied to the path $v_3P_3w_3$
there exist two disjoint paths $Q_1,Q_2$ in $G$ from $Y$
to $V(P_1\cup P_2)\cup\{w_3\}$,
with no internal vertices in $Y\cup V(P_1\cup P_2)$ and
such that one of them, say $Q_2$, ends in $w_3$.
From the symmetry we may assume that $Q_1$ ends in $V(P_1)$.
Similarly as before, let $C_{12}$ denote the subpath of $C\backslash u_3$
with ends $u_1$ and $u_2$, and let $C_{13}$ and $C_{23}$ be defined
similarly.
If $C_{23}$ has at least $\ell/5$ edges, then the disjoint subgraphs
$P_1\cup Q_1$ and $P_2\cup Q_2\cup C_{23}\cup w_3P_3u_3$ include
a linkage from $X$ to $Y$ of length at least $\ell/5$, as desired.
Thus we may assume that $|E(C_{23})|<\ell/5$.
If $u_1\not\in V(Q_1)$, then replacing $C_{23}$ by $C_{12}\cup C_{13}$
above results in a linkage from $X$ to $Y$ of length at least $4\ell/5$.
Thus we may assume that $V(Q_1)\cap V(P_1)=\{w_1\}$ and $w_1=u_1$.
But now either
$(P_1\cup Q_1,P_2\cup C_{23}\cup Q_2,C_{12},C_{13})$
is a hammock from $X$ to $Y$ of length at least $2\ell/5$, or
$(P_1\cup Q_1,P_2\cup C_{23}\cup Q_2,C_{13},C_{12})$
is a hammock from $Y$ to $X$ of length at least $2\ell/5$,
as desired.~\qed

\begin{lemma}
\label{lem:C}
Let $G$ be a $3$-connected graph, let $X,Y\subseteq V(G)$ be
disjoint sets of size two, and let $G$ have a hammock from
$X$ to $Y$ of length $l$.
Then $G$ has either a non-singular hammock from $X$ to $Y$
or from $Y$ to $X$ of length at least $l/2$, or
a linkage from $X$ to $Y$ of length at least $l/2$.
\end{lemma}

\paragraph{Proof:}
Let $\eta=(P_1,P_2,R_1,R_2)$ be a hammock in $G$ of length
$l$, and let $x_1,x_2,y_1,y_2,s_1,s_2,t_1,t_2$
be as in the definition of hammock.
We may assume that $s_1=s_2$, for otherwise $\eta$ is
non-singular, and hence satisfies the conclusion of the lemma.
Since $x_1\ne y_1$, at least one of the sets
$A:=V(x_1P_1s_1)-\{s_1\}$ and $B:=V(y_1P_1s_2)-\{s_2\}$
is not empty.

Assume first that $A\ne\emptyset$.
Since $G$ is $3$-connected, there is a path $Q$ in
$G\backslash\{s_1,t_1\}$ with ends $a\in A$ and
$b\in V(P_2\cup R_1\cup R_2\cup s_2P_1y_1)$.
If $b\in V(R_2\cup t_1P_2y_2)$, then
$x_1P_1s_1\cup Q\cup R_2\cup t_1P_2y_2$ includes a path from
$x_1$ to $y_2$ that together with
$x_2P_2t_1\cup R_1 \cup s_2P_1y_1$ forms a linkage from
$X$ to $Y$ of length at least $l$.
If $b\in V(x_2P_2t_1)$, then $(P_1,P_2,R_1,Q)$ is a
non-singular hammock
from $Y$ to $X$ of length $l$.
If $b\in V(s_1P_1y_1)$, then the paths $x_1P_1a\cup Q\cup bP_1y_1$
and $x_2P_2t_1\cup R_1\cup R_2\cup t_2P_2y_2$ form a linkage
from $X$ to $Y$ of length at least $l$.
Thus we may assume that $b\in V(R_1)$. If the path $t_1R_1b$
has at least $l/2$ edges, then
$(P_1,P_2,Q\cup t_1R_1b,R_2)$ is a non-singular hammock from $X$ to $Y$
of length at least $l/2$,
and if the path $s_1R_1b$ has at least $l/2$ edges, then
the paths $P_2$ and $x_1P_1a\cup Q\cup bR_1s_1\cup s_1P_1y_1$
form a linkage from $X$ to $Y$ of length at least $l/2$.
This completes the case $A\ne\emptyset$.

Thus we may assume that $B\ne\emptyset$.
%If $B\ne\emptyset$ then
We take a path in
$G \setminus \{s_2,t_2\}$ connecting a vertex in $B$ to a vertex in
$V(P_2 \cup R_1 \cup R_2 \cup x_1P_1s_1)$ and
proceed similarly as in the previous paragraph.
%then consider the different locations of the end point
%in $V(P_2 \cup R_1 \cup R_2 \cup x_1P_1s_1)$ as above.
The details are analogous to the case $A\ne\emptyset$
and are left to the reader. ~\qed
\bigskip

\begin{lemma}
\label{lem:B}
Let $G$ be a $3$-connected graph, let $X,Y\subseteq V(G)$ be
disjoint sets of size two, and assume that $G$ has a non-singular
hammock from $X$ to $Y$ of length $l$.
Then $G$ has a linkage from $X$ to $Y$ of length at least $l/2$.
\end{lemma}

\paragraph{Proof:}
Let $\eta=(P_1,P_2,R_1,R_2)$ be a non-singular hammock in $G$,
and let $x_1,x_2,y_1,y_2,s_1,s_2,t_1,t_2$
be as in the definition of hammock.
Since $G$ is non-singular, $(P_2,P_1,R_1,R_2)$ is also a hammock
from $X$ to $Y$ of the same length as $\eta$, and hence there is
symmetry between $P_1$ and $P_2$.
Let $A:=x_1P_1s_1\cup R_1\cup x_2P_2t_1$ and
$B:=y_1P_1s_2\cup R_2\cup y_2P_2t_2$.
Then $A$ and $B$ are paths in $G$.
By Lemma~\ref{augment} applied to the sets $V(A)$ and $V(B)$ and paths
$s_1P_1s_2$ and $t_1P_2t_2$ there exist three disjoint paths
$Q,Q_1,Q_2$ from $V(A)$ to $V(B)$ such that two of them have ends
in $\{s_1,t_1\}$, and two have ends in $\{s_2,t_2\}$.
From the symmetry we may assume that $s_1$ is an end of $Q_1$;
let $s_2'$ be the other end of $Q_1$.
Similarly we may assume  that $t_1$ is an end of $Q_2$;
let $t_2'$ be the other end of $Q_2$.
Now the path $x_1P_1s_1 \cup Q_1\cup s_2'By_1$ can play the role of $P_1$,
the path $x_2P_2t_1\cup Q_2\cup t_2'By_2$ can play the role of $P_2$,
and the path $s_2'Bt_2'$ can play the role of $R_2$.
In other words, we may assume (by changing the hammock $\eta$ but not
changing its length) that there exists a path $Q$ from
$V(A)$ to $V(B)$ that is disjoint from the paths $P_1$ and $P_2$.
Let $a$ be the end of $Q$ in $V(A)$, and let $b$ be the end of $Q$ in $V(B)$.
From the symmetry we may assume that either $a\in V(x_1P_1s_1)$,
or $a\in V(R_1)$ and the path $aR_1t_1$ has at least $l/2$ edges.

Assume first that $a\in V(x_1P_1s_1)$.
If $b\in V(s_2P_1y_1)$, then the paths $x_1P_1a\cup Q\cup bP_1y_1$ and
$x_2P_2t_1\cup R_1\cup s_1P_1s_2\cup R_2\cup t_2P_2y_2$ form a linkage
from $X$ to $Y$ of length at least $l$, and if $b\in V(R_2\cup t_2P_2y_2)$,
then the path $x_2P_2t_1\cup R_1\cup s_1P_1y_1$ and a subpath of
$x_1P_1a\cup Q\cup R_2\cup t_2P_2y_2$ form a linkage
from $X$ to $Y$ of length at least $l$.
This completes the case when $a\in V(x_1P_1s_1)$.

We may therefore assume that $a\in V(R_1)$ and
the path $aR_1t_1$ has at least $l/2$ edges.
If $b\in V(s_2P_1y_1)$, then the paths $x_1P_1s_2\cup R_2\cup t_2P_2y_2$
and $x_2P_2t_1\cup t_1R_1a\cup Q\cup bP_1y_1$ form a linkage
from $X$ to $Y$ of length at least $l/2$, and if $b\in V(R_2\cup t_2P_2y_2)$,
then the path $P_1$ and a subpath of
$x_2P_2t_1\cup t_1R_1a\cup Q\cup R_2\cup t_2P_2y_2$ form a linkage
from $X$ to $Y$ of length at least $l/2$.~\qed

\paragraph{Proof of Lemma~\ref{disjointpaths}.}
Let $G,X,Y$ be as stated, and assume that $G$ has a path of length $\ell$.
Then $G$ has a cycle of length at least $2\ell/5$ by Theorem~\ref{BLtheo}.
By Lemma~\ref{lem:A} we may assume that $G$ has a hammock from $X$
to $Y$ of length at least $4\ell/25$, for otherwise the theorem holds.
Similarly, by Lemma~\ref{lem:C} we may assume that $G$ has a non-singular
hammock from $X$ to $Y$ of length at last $2\ell/25$.
By Lemma~\ref{lem:B} the graph $G$ has a linkage from $X$ to $Y$
of length at least $\ell/25$, as desired.~\qed

\section{Proof of Theorem~\ref{theomain}}
%%%%%%%%%%%%%%%%%%%%%%%%%%%%%%
\label{secmain}

%\begin{claim}\label{logs}
%For every $y \geq x \geq 2$ we have $\log x + \log y \geq \log (x+y)$
%\end{claim}

%\paragraph{Proof:} Just observe that $\log (x+y)-\log y \leq \log 2y - \log y = \log 2 \leq \log x$. $\qed$

\paragraph{Notation.} Throughout this section we will assume the
following notation.
Let $k\ge4$ be an integer, let $G$ be a $k$-critical graph
on $n$ vertices, and
let $(T,{\cal W})$ be a standard tree-decomposition of $G$.
One exists by Lemma~\ref{standard}.
For $t\in V(T)$ let $H_t$ denote the torso of $(G,T,{\cal W})$ at $t$,
and let $N_t$ denote the nucleus of $(G,T,{\cal W})$ at $t$.
We select a vertex $r\in V(T)$ of degree one that we will regard
as the root of $T$. Thus  a {\em descendant} of a vertex $t\in V(T)$ is
any vertex $t'\in V(T)-\{t\}$ such that $t$ belongs to the path
from $r$ to $t'$ in $T$.
For $t\in V(T)$ we denote by $T_t$ the subtree of $T$ induced by $t$ and
all its descendants.
We define a weight function $w:V(T)\to\{0,1,\ldots\}$
by $w(t):=|E(N_t)|$.
Thus $w(t)\ge 6$ for every $t\in V(T)$ by Lemma~\ref{nucleus}(ii).
According to the convention introduced prior to lemma~\ref{longpath},
$w(T_t)$ means $\sum_{v\in V(T_t)}w(t)$.
We now define, for every $t\in V(T)$, a set $X_t\subseteq W_t$ of size two.
If $t \neq r$, then let $t'$ be the parent of $t$ in the rooted tree $(T,r)$
and we set $X_t=W_t \cap W_{t'}$.
If $t=r$ and $T$ has at least two vertices, then let $t'$ be the unique
child of $r$ in $T$ and let
$X_r\subseteq W_r$ be any set disjoint from $W_t\cap W_{t'}$
that consists of two vertices that are adjacent in $G$.
Such a set exists because $H_t$ is $3$-connected by Lemma~\ref{3connbags}
and the elements of $W_t\cap W_{t'}$ form the only edge of $H_t$ that
does not belong to $G$.
Finally, if $T$ has only one vertex we choose $X_t$ arbitrarily.
For $t\in V(T)$ we denote by $G_t$ the graph induced in $G$ by the set
of vertices $\bigcup_{t' \in V(T_t)}W_{t'}$.

In order to be able to apply Lemma~\ref{longpath} we prove the
following lemma.

\begin{lemma}
\label{sptree}
Let $t\in V(T)-\{r\}$, and let $X_t=\{x,x'\}$.
Then the graph $G_t \setminus x$ has a spanning tree $R$
such that for every integer $\ell\ge0$ there are at most $k^\ell$
vertices of $R$ at distance exactly $\ell$ from $x'$.
\end{lemma}

\paragraph{Proof.}
%To complete the induction proof in case (\ref{secondcase}) holds, we
%need to show the existence of a spanning tree $R$ of $G_t \setminus x$
%as described above. First, observe that the weight function $\phi$ we
%have defined above satisfies $\phi(x')=0$ and $\phi(v) \neq 0$ for at
%least one vertex so we just need to show that $G_t \setminus x$ has a
%spanning tree rooted at $x'$ with at most $k^{j+1}$ vertices at depth $j$.
%As before, let $G_t=G[\bigcup_{v \in T_t}W_t]$ and
Let $G'_{t}$ be the subgraph of $G$ induced by the union of all
$W_{t'}$ over all $t' \in V(T)-V(T_t)$.
Then $G_{t} \cap G'_{t} = X_t=W_t \cap W_{t'}$, where $t'$ is the
ancestor of $t$ in the rooted tree $(T,r)$.
By Lemma \ref{2con}(iv)
applied to $G$ and the vertices $x,x'$
the graph $G_{t}$ was obtained from some $k$-critical graph $H$ by either
(i) deleting the edge $xx'$, or
(ii)  splitting  a vertex of $H$ into the two vertices $x,x'$.

Assume first that $G_{t}$ was obtained by deleting the edge $xx'$.
Since $H$ is $2$-connected by Lemma~\ref{2con}(i), the graph $H \setminus x$
has a DFS spanning tree $R$ rooted at $x'$. We deduce
from Lemma~\ref{AKSmod} applied to $H$ and $X=\{x\}$ that $R$ satisfies
the conclusion of the lemma.

We may therefore assume that $G_{t}$ was obtained from $H$ by splitting a
vertex, say  $z$, into the two vertices $x,x'$.
Since $G$ is $2$-connected by Lemma~\ref{2con}(i), there is an edge
$e\in E(G_t)$ joining the vertex $x'$ to a vertex in $V(G_t)-\{x,x'\}$.
Then $e$ is also an edge of $H$.
Since $H$ is $2$-connected, it has a DFS spanning tree $R'$ rooted at $z$
such that $e$ is the only edge of $R'$ incident with $z$.
The tree $R'$ gives rise to a unique spanning tree $R$ of $G_t\backslash x$
with the same edge-set in the obvious way. It follows
from Lemma~\ref{AKSmod} applied to $H$ and $X=\emptyset$ that $R$ satisfies
the conclusion of the lemma.~\qed
\bigskip

To prove Theorem~\ref{theomain} we prove, for the sake of induction,
 the following lemma.

\begin{lemma}
\label{keylemma}
For every $t \in V(T)$ the graph $G_t$ has a path connecting
the vertices of $X_t$ of length at least $\log w(T_t)/(100\log k)$.
%$\frac{1}{36C}\log w(T_t)/\log k$.
\end{lemma}

Let us first derive Theorem~\ref{theomain} from Lemma~\ref{keylemma}.

\paragraph{Proof of Theorem \ref{theomain}, assuming Lemma~\ref{keylemma}.}
Apply Lemma \ref{keylemma} with $t=r$. Since by definition
the vertices of $X_r$ are adjacent,
we get a cycle of length at least $\log w(T)/(100\log k)$.
Since distinct nuclei are edge-disjoint by the definition of nucleus,
Lemma~\ref{edgesinnuc} implies
$w(T)=\sum_{t\in V(T)}|E(N_t)|\ge |E(G)|\ge n$,
and hence $G$ has a cycle of
length at least $\log n/(100\log k)$.~$\qed$

\bigskip

The rest of this section is devoted to a proof of Lemma \ref{keylemma}.
We first take care of the following special case.

\begin{lemma}
\label{basecase} Let $t\in V(T)$.
The statement of Lemma~\ref{keylemma} holds for $t$ if
$w(t) \geq w(T_t)/5$.
In particular, the lemma holds for $t$ if $|V(T_t)|=1$.
\end{lemma}

% and observe that since $t$ is a leaf of $T$ we infer that $H_t$ is obtained from $G[W_t]$ be adding an edge connecting the vertices of $X_t$

\paragraph{Proof.}
The second assertion follows from the first, and so it suffices to
prove the first statement.
By Lemma~\ref{pathtorso} the torso $H_t$ has a path of length at least
$\frac{1}{2}\log w(t)/\log k$ (recall
that $w(t)$ is the number of edges in the nucleus  $N_t$).
Lemma \ref{3connbags} guarantees that $H_t$ is $3$-connected,
and so by Theorem \ref{BLtheo}
we get that $H_t$ has a cycle $C$ of length at least
$\frac{1}{5}\log w(t)/\log k$.
Since $H_t$ is $3$-connected we get from Menger's Theorem that it
contains two disjoint paths connecting $X_t$ to $C$. Suppose these paths meet
$C$ at vertices $w,w'$. Then one of the two subpaths of $C$ connecting
$w$ and $w'$ has length at least
$\frac{1}{10}\log w(t)/\log k$.
Together with the two paths connecting the vertices of $X_t$ to $C$
we get a path $P$ in $H_t$ connecting the vertices
of $X_t$ of length at least
$\frac{1}{10}\log w(t)/\log k\geq \frac{1}{100}\log w(T_t)/\log k$
by the hypothesis of the lemma and the fact that $w(T_t)\ge6$.
%because every nucleus has at least six edges by Lemma~\ref{nucleus}(ii).
%Since the lemma assumes that $w(t) \geq \frac12 w(T_t)$ we get a path of length at least $\frac{1}{36C}\log w(T_t)/\log k$.

For every edge $e=uv\in E(P)-E(G)$ we do the following.
By Lemma~\ref{virtue}(ii) there is a unique neighbor $t'$ of $t$ in $T$
such that $W_t \cap W_{t'}=\{u,v\}$.
If $r\ne t$, then $t'$ is not the parent of $t$ in the rooted tree
$(T,r)$, because $\{u,v\}\ne X_t$ by the choice of $P$.
%The graph $G$ is $2$-connected by Lemma~\ref{2con}(i).
We claim that there exists a path $P_e$ in $G_{t'}$ with ends $u,v$.
Indeed, since $(T,{\cal W})$ is standard, there exists a vertex
$w\in W_{t'}-\{u,v\}$.
Since $G$ is $2$-connected by Lemma~\ref{2con}(i), there exist two
paths $P_1,P_2$ in $G$ with one end $w$ and the other end in $\{u,v\}$,
pairwise disjoint, except for $w$.
Then $P_e:=P_1\cup P_2$ is a path in $G$ with ends $u,v$.
It follows that $P_e$ is a path in $G_{t'}$, for otherwise some subpath $Q$
of $P_e\backslash\{u,v\}$ joins the vertex  $w\in V(G_{t'})-\{u,v\}$
to a vertex of $V(G)-V(G_{t'})$.
But then $Q$ has an edge with one end in $V(G_{t'})-\{u,v\}$ and the other
end in $V(G)-V(G_{t'})$, contrary to definition of tree-decomposition.
This proves our claim that $P_e$ exists.
We replace $e$ by the path $P_e$ and repeat the construction for
each edge $e\in E(P)-E(G)$.
For distinct edges $e,e'\in E(P)-E(G)$ the paths $P_e,P_{e'}$ have no
internal vertices in common, because their interiors belong to
disjoint subgraphs.
%by Lemma~\ref{virtue}(ii).
We thus arrive at a path in $G$ with ends
in $X_t$ of length at least $\log n/(100\log k)$, as desired.~$\qed$
\bigskip

\paragraph{Proof of Lemma~\ref{keylemma}.}
%We complete the proof of Lemma \ref{keylemma} \
We proceed by induction on $|V(T_t)|$.
Since Lemma~\ref{basecase} establishes the base case $|V(T_t)|=1$,
we can assume henceforth that $|V(T_t)|>1$ and that the lemma holds for
trees of size less than $|V(T_t)|$.
Let $X_t=\{x,x'\}$, let $N$ be the children of $t$ in the rooted tree
$(T,r)$ and define
$$
N_0=\{t' \in N~:~W_{t}\cap W_{t'}\cap X_t= \emptyset\}
$$
$$
N_1=\{t' \in N~:~W_{t}\cap W_{t'}\cap X_t= \{x\}\}
$$
$$
N_2=\{t' \in N~:~W_{t}\cap W_{t'}\cap X_t= \{x'\}\}
$$
The sets $N_0,N_1,N_2$ form a partition of $N$. For $t\ne r$ this follows
from Lemma~\ref{virtue}(ii), and for $t=r$
this follows from the way we picked $X_r$.
Therefore, either
%some $i \in \{0,1,2\}$ we have
\begin{equation}
\label{eqcase1}
\sum_{y \in N_0}w(T_{y}) \geq \frac34(w(T_t)-w(t)),
\end{equation}
or
\begin{equation}
\label{eqcase2}
\sum_{y \in N_1\cup N_2}w(T_{y}) \geq \frac14(w(T_t)-w(t)).
\end{equation}
We first deal with the case (\ref{eqcase1}).
By Lemma~\ref{basecase} we may assume that $w(t) <  w(T_t)/5$, and hence
%already handles the case in which $w(t) \geq \frac12 w(T_t)$ we can actually assume that
%$$
%\sum_{y \in N_i}w(T_{y}) \geq \frac16 w(T_t).
%$$
%From the symmetry between the cases $i=1$ and $i=2$ it is enough to
%consider only the cases $i=0$ and $i=1$. Let us first consider
%the case $i=0$; that is
\begin{equation}\label{firstcase}
\sum_{y \in N_0}w(T_{y}) \geq \frac35 w(T_t)\;.
\end{equation}
By Lemma~\ref{degreeT} we know that $|N_0| \leq |N| \leq 3w(t)$.
Therefore, there is a vertex $t' \in N_0$ for which
\begin{equation}\label{t*}
w(T_{t'})\geq \frac{w(T_t)}{5w(t)}\;.
\end{equation}

\noindent
By Lemma~\ref{pathtorso} the graph $H_t$ (the torso at $t$) has a path of
length at least $\frac{1}{2}\log w(t)/\log k$. Therefore,
by Lemma~\ref{3connbags}, we can apply
Lemma~\ref{disjointpaths} to the graph $H_t$ and sets $X_t$ and
$X_{t'}$ to deduce that $H_t$ has two disjoint paths $P_1,P_2$ from
$X_t$ to $X_{t'}$ satisfying
\begin{equation}\label{P1P2}
|E(P_1)|+|E(P_2)| \geq  \frac{\log w(t)}{50 \log k}\;.
\end{equation}
By the induction hypothesis  the graph $G_{t'}$
%$G[\bigcup_{v \in T_{y^*}}W_v]$
has a path $P$ connecting the pair of vertices of $X_{t'}$ satisfying
\begin{equation}\label{P}
|E(P)| \geq \frac{\log w(T_{t'})}{100\log k}\;.
\end{equation}
Combining (\ref{t*}), (\ref{P1P2}) and (\ref{P}) we get that
$P_1 \cup P \cup P_2$ is a path in $G_t\cup H_t$ with ends in $X_t$ of length
at least
\begin{eqnarray*}
|E(P_1)|+|E(P)|+|E(P_2)| &\geq&
\frac{\log w(t)}{50\log k}+\frac{\log w(T_{t'})}{100\log k}\\
%&\geq& \frac{\log w(t)}{50\log k}+\frac{\log ({w(T_t)/18w(t)})}{100\log k}\\
&\geq& \frac{\log w(t)}{50\log k}+\frac{\log {w(T_t)}}{100\log k}-
\frac{\log {w(t)}+\log5}{100\log k}\\
&\geq& \frac{\log w(T_t)}{100\log k}\;,
\end{eqnarray*}
because $w(t)\ge6$.
We now convert $P_1\cup P\cup P_2$ to a path in $G_t$ of length
at least $\log{w(T_t)}/(100\log k)$ in the same
way as in the second paragraph of the proof of Lemma~\ref{basecase}.
This completes the proof when (\ref{eqcase1}) holds.

Thus we may assume~(\ref{eqcase2}).
From the symmetry between $N_1$ and $N_2$ we may assume that
$$%\begin{equation}
%\label{eqcase2}
\sum_{y \in N_1}w(T_{y}) \geq (w(T_t)-w(t))/8.
$$%\end{equation}
Again, by Lemma~\ref{basecase}
we may assume that $w(t) <  w(T_t)/5$, and hence
\begin{equation}\label{secondcase}
\sum_{y \in N_1}w(T_{y}) \geq w(T_t)/10\;.
\end{equation}

%We now return to the proof of the induction when (\ref{secondcase}) holds.
It follows that $t \neq r$, for otherwise $N_1 = \emptyset$.
We need to define a new weight function $\phi:V(G_t) -\{ x\}\to\{0,1,\ldots\}$.
Let $v\in V(G_t) -\{ x\}$.
If $v\in W_t$ and there exists a neighbor $t'$ of $t$ in $T_t$ such that
$W_{t'} \cap W_t=\{x,v\}$, then $t'$ is unique by Lemma~\ref{virtue}(ii),
and we define $\phi(v)=w(T_{t'})$. If $v\not\in W_t$ or no
such $t'$ exists, then we define $\phi(v)=0$.
Thus, in particular, $\phi(x')=0$ by Lemma~\ref{virtue}(ii).
By Lemma~\ref{sptree} the graph $G_t \setminus x$ has a spanning tree $R$
such that for every integer $l\ge0$ there are at most $k^\ell$ vertices
of $R$ at distance exactly $\ell$ from $x'$.
Note that by (\ref{secondcase}) we have that the total weight of $R$ satisfies
\begin{equation}\label{weightR}
\phi(R)= \sum_{y \in N_1} w(T_{y})\geq  w(T_t)/10\;.
\end{equation}
By Lemma~\ref{longpath} applied to the tree $R$ and vertex $x'$
there exists a vertex $v\in V(R)$ at distance $\ell$ from $x'$ in $R$
such that $\phi(v)>0$ and $2\ell\log k+\log\phi(v)\ge\log\phi(T)$.
It follows that there is a path $P$ in $G_t \setminus x$
 from $x'$ to $v$ satisfying
\begin{equation}\label{boundP}
2|E(P)|\log k + \log \phi(v) \geq \log \phi(R)\;.
\end{equation}
%Note that since only vertices in $W_t \setminus \{x,x'\}$ have non-zero weight it must be one of these vertices.
%Combining (\ref{weightR}) and the fact that $\phi(v)=w(T_{t'})$ we get that there is a path $P$ from $x'$ to $v$ that avoids $x$ and satisfies
Since $\phi(v)>0$ we deduce that $v\in W_t$ and $P$ has length at least one.
Let $t'\in N_1$ be such that $W_{t'} \cap W_t=\{x,v\}$, so that
$\phi(v)=w(T_{t'})$.
Since $P$ is a path from $x'\in W_t-V(G_{t'})$ to $v$ in $G_t\backslash x$,
we deduce that $V(P)\cap V(G_{t'})=\{v\}$.
By the induction hypothesis applied to the graph $G_{t'}$ the graph
 $G_{t'}$ has a path $Q$ connecting $x$ to
$v$ of length at least $\frac{1}{100}\log \phi(v)/\log k$.
So $P \cup Q$ is a path in $G_t$ from $x$ to $x'$ of length at least
\begin{eqnarray*}
|E(P)|+\frac{\log \phi(v)}{100\log k} &=&
\frac{1}{50}\left(|E(P)|+\frac{\log \phi(v)}{2\log k}\right)+
\left(1-\frac{1}{50}\right)|E(P)|\\
&\geq& \frac{\log \phi(R)}{100\log k}+1-\frac{1}{50} \\
&\geq& \frac{\log (w(T_t)/10)}{100\log k}+1-\frac{1}{50}\\
&=& \frac{\log w(T_t)}{100\log k}-\frac{\log 10}{100\log k}+1-\frac{1}{50}\\
&\geq& \frac{\log w(T_t)}{100\log k}\;,
\end{eqnarray*}
where in the first inequality we used (\ref{boundP})
and the fact that $|E(P)| \geq 1$,
and the second inequality uses (\ref{weightR}).
This completes the proof of Lemma~\ref{keylemma}.~\qed

\section{Gallai's Upper Bound}\label{secgallai}

We need to introduce the notion of Haj\'os sum of two graphs.
Let $K$ and $L$ be two graphs with disjoint vertex-sets, and let
$k_1k_2$ and $l_1l_2$ be edges of $K$ and $L$, respectively.
Let $G$ be the graph obtained from the union of $K$ and $L$ by
deleting the edges $k_1k_2$ and $l_1l_2$, identifying the vertices
$k_1$ and $l_1$, and adding an edge joining $k_2$ and $l_2$.
In those circumstances we say that $G$ is a {\em Haj\'os sum} of $K$ and $L$.
It is straightforward to check that if $K$ and $L$ are $k$-critical,
then so is $G$.

We now describe a construction of $k$-critical graphs with no
long path, and hence no long cycle.
Let $k\ge4$ be an integer, let $T$ be a tree of maximum degree at most
$k-1$, and let $(H_t:t\in V(T))$ be a family of $k$-critical graphs,
each containing the same vertex $x_0$, and otherwise pairwise disjoint.
For every ordered pair $t,t'$ of adjacent vertices of $T$ we select
a vertex $v_{tt'}\in V(H_t)$ such that
\myitem{$\bullet$} $v_{tt'}$ is adjacent to $x_0$ in $H_t$, and
\myitem{$\bullet$} if $t'$ and $t''$ are distinct neighbors of $t$ in $T$, then
$v_{tt'}\ne v_{tt''}$.

\noindent
Such a choice is possible, because $T$ has maximum degree at most $k-1$
and every $k$-critical graph has minimum degree at least $k-1$.
Let us emphasize that even though $tt'$ and $t't$ denote the same
edge of $T$, the vertices $v_{tt'}$ and $v_{t't}$ are distinct:
the first belongs to $H_t$ and the second to $H_{t'}$.
We define a graph $G$ to be the graph obtained from $\bigcup_{t\in V(T)} H_t$
by, for every edge $tt'\in E(T)$, deleting the edges
$x_0v_{tt'}$ and $x_0v_{t't}$,
and adding an edge joining $v_{tt'}$ and $v_{t't}$.

% deleting all the edges $x_0v_{tt'}$ for all ordered pairs
%of adjacent vertices $t,t'\in V(T)$, and adding, for every edge $tt'\in E(T)$,
%an edge joining $v_{tt'}$ and $v_{t't}$.

It is easy to see that the graph $G$ can be viewed as being obtained
from the graphs $H_t$ by repeatedly taking Haj\'os sums, and is thus
$k$-critical. Also, it has $1+\sum_{t\in V(T)}(|V(H_t)|-1)$ vertices, and for
every path $P$ in $G\backslash x_0$ there exists a path $R$ in $T$ such that
$|V(P)|\le \sum_{t\in V(R)}(|V(H_t)|-1)$.
To replicate (a close relative of) Gallai's original construction, let
$h\ge0$ be an integer, and let $T$ be the $(k-1)$-branching tree of height $h$;
that is, a tree $T$ with a vertex $r\in V(T)$ such that every vertex
is at distance at most $h$ from $r$, and each vertex at distance at most
$h-1$ from $r$ has degree exactly $k-1$. Each of the graphs $H_t$ will be
the complete graph on $k$ vertices.
Then the graphs resulting from the construction described above with
this choice of $T$ and $H_t$ prove the inequality~(\ref{upperbound})
and statement~(\ref{gallaipath}),
as is easily seen.
(In Gallai's original construction the vertex $r$ has degree $k-2$,
but that makes little difference.)
However, there exist $k$-critical graphs on $n$ vertices for every
$n\ge k$, except $n=k+1$.
It is easy to deduce the following theorem,
by utilizing such $k$-critical graphs and trees that are not
necessarily regular, and the above construction.

\begin{theo}
\label{supergallai}
For every integer $k\ge4$ and every integer $n\ge k+2$
\begin{equation}
\label{newgallai}
L_k(n) \leq \frac{2(k-1)}{\log (k-1)}\log n+2k\;.
\end{equation}
\end{theo}

\bigskip

\paragraph{Acknowledgment.}
We are grateful to Michael Krivelevich for answering many questions related
to this paper and especially for sharing
his English translation of the construction of Gallai from \cite{G}.

\baselineskip 11pt
\vfill
\noindent
This material is based upon work supported by the National Science Foundation.
Any opinions, findings, and conclusions or
recommendations expressed in this material are those of the authors and do
not necessarily reflect the views of the National Science Foundation.
\eject

\end{document}